\documentclass[12pt]{article}
\usepackage{fullpage, amssymb, latexsym, amsmath, graphicx,
epsfig, amsfonts, amsthm, bm, float}
\usepackage[numbers]{natbib}

\newtheorem{theorem}{Theorem}
\newtheorem{corollary}{Corollary}
\newcommand{\bR}{\mathbb{R}}
\newcommand{\bZ}{\mathbb{Z}}
\DeclareMathOperator{\E}{E}

\newenvironment{problist}
   {
      \begin{list}
         {---}
         {
            \setlength{\itemsep}{.5ex}
            \setlength{\parsep}{0ex}
            \setlength{\parskip}{0ex}
            \setlength{\topsep}{.5ex}
         }
   }
   {
      \end{list}
   }

\begin{document}

\bibliographystyle{plainnat}

\date{}

\title{Multivariate Lag-Windows and Group Representations\thanks{Last Updated: \today}}

\author{{\Large Arthur Berg}\\{\normalsize {\it University of California, San Diego}}\\
\texttt{aberg@ucsd.edu} }

\maketitle \abstract{ Symmetries of the auto-cumulant function (the
generalization of the auto-covariance function) of a
$k^{\text{th}}$-order stationary time series are derived through a
connection with the symmetric group of degree $k$. Using theory of
group representations, symmetries of the auto-cumulant function are
demystified and lag-window functions are symmetrized to satisfy
these symmetries. A generalized Gabr-Rao optimal kernel, used to
estimate general $k^{\text{th}}$-order spectra, is also derived
through the developed theory.
 }\vspace{.5cm}


\section{Introduction}
There have been numerous approaches and references to constructing
multivariate lag-windows that inherit the appropriate symmetries
when estimating $k^{\text{th}}$-order spectra; cf. \cite{Lu96,
Nikias87, Nikias93b,Nikias93,Sasaki75,rao84, rao03,
Turkbeyler00,Yao01}. However, each approach constructs multivariate
lag-windows from univariate lag-windows which severely limits their
shape.  The case $k=2$, corresponding to univariate lag-windows, is
trivial since the only symmetry requirement imposed is evenness. The
symmetries are much more complex in the case $k=3$, and the
challenge is to fully understand these symmetries and to easily
construct lag-windows that inherit the symmetries. In the process,
an intimate relationship is discovered between the symmetries of the
bivariate lag-window and the symmetric group $S_3$.  Here we begin
to see the makings of a beautiful connection between stationary time
series and group representations of the symmetric group.

General theory of the symmetric group is pulled back via this
connection to the context of stationary time series and many useful
results are ascertained.  Similarly, natural constructions of
stationary time series are pushed forward into the theory of group
representations of $S_k$ to elicit a new realization of a familiar
representation of $S_k$ with dimension $k-1$.

Current practices of generating multivariate lag-windows are shown
to be inadequate as they restrict the construction of many useful
lag-windows including the Gabr-Rao optimal bivariate lag-window
\cite{rao84}, yet this lag-window occurs naturally under the
symmetrization routine described in this paper.  The connection to
the symmetric group via group representations allows for effortless
generalization of the Gabr-Rao lag-window for estimation of the
$k^{\text{th}}$-order spectral density.

This is not the first time a connection between group
representations and statistics has been drawn.  In fact an entire
book {\it Group Representations in Probability and Statistics}
\cite{Diaconis88}, written by the mathemagician Persi Diaconis,
exposes many connections between group representations and
statistics, but the connections drawn in this book are much
different, with different motivations, from what is done the present
article.

The first section defines stationarity---the driving force in the
connection with the symmetric group---and introduces higher-order
spectral densities and spectral density estimates with lag-windows.
To solidify ideas, the case $k=3$ is considered in Section 2, and a
connection between symmetries of the auto-cumulant function (acf)
and permutations of the symmetric group is obtained via a special
group representation on $S_3$; these results are generalized for
general $k$ in Section 4. Section 5 describes the actions in the
Fourier domain viewpoint which allows for the generalization of the
Gabr-Rao kernel. Finally, several examples of bivariate kernels and
lag-windows are produced in Section 6.

\section{Stationarity}

Let $\{X_t\}$ be a $k^{\text{th}}$-order discrete (resp. continuous)
stationary time series, that is
\begin{equation}
\label{eqn:expectation} \E\left[X_t\prod_{t_i\in S}X_{t+x_i}\right]
\end{equation}
is finite and the same for every $t\in\bZ$ (resp. $t\in\bR$) where
$S\subset\bZ$ (resp. $S\subset\bR$) is any set of size at most
$k-1$.  A stronger, yet more typical, assumption is to require the
time series to be strictly stationary, that is for any positive
integer $n$ and any $t, x_1,\ldots,x_n$ in $\bZ$ (resp. $\bR$), the
random vectors $(X_{x_1},\ldots,X_{x_n})$ and
$(X_{t+x_1},\ldots,X_{t+x_n})$ have identical joint distributions.
However, this stronger version of stationarity will not be needed
here since this paper only deals with expressions like those in
(\ref{eqn:expectation}).

As far as this paper is concerned, the theory for discrete time
series is the same as that for continuous time series, so the
remainder of this paper will focus on discrete time series only.

Let's specialize for the moment to the case of second-order (aka
weakly) stationary time series.  The assumption of second-order
stationarity requires
\begin{problist}
  \item[(i)]  $\E(X_t^2)<\infty$ for all $t$;
  \item[(ii)] $\E(X_t)=\mu$ for all $t$;
  \item[(iii)] $C(x)=\E[(X_t-\mu)(X_{t+x}-\mu)]$ for all $x$ and $t$.
\end{problist}
The function $C(x)$ is called the auto-covariance function (referred
to within as the auto-cumulant function), and by replacing $x$ with
$-x$, we see that
\[
C(-x)=\E[(X_t-\mu)(X_{t-x}-\mu)]=\E[(X_{t'}-\mu)(X_{t'+x}-\mu)]=C(x)
\]
(using $t'=t+x$); so $C(x)$ is an even function, and this is the
only symmetry requirement of $C(x)$.

Third-order stationarity is implied by the conditions (i)-(iii)
above plus following two conditions:
\begin{problist}
  \item[(iv)] $\E(X_t^4)<\infty$
  \item[(v)] $C(x,y)=\E[(X_t-\mu)(X_{t+x}-\mu)(X_{t+y}-\mu)]$ for all $x$, $y$, and $t$.
\end{problist}
In exploring the symmetries of the bivariate acf, $C(x,y)$, we can
use the substitutions $t'=t-x$ and $t'=t-y$ and the commutativity
property of multiplication to deduce
\begin{equation}
\label{eqn:bivariatesym}
  C(x,y)=C(-x,y-x)=C(y,x)=C(x-y,-y)=C(-y,x-y)=C(y-x,-x)
\end{equation}
There is redundancy in this expression; for instance the relations
$C(x,y)=C(y,x)=C(-x,y-x)$ imply $C(x,y)=C(y-x,-x)$.  This begs two
natural questions:
\begin{problist}
  \item  Which equations can be removed without losing any
  information?
  \item  Is there minimal set of equations representing all of the
  symmetries?
\end{problist}
We will answer both of these questions in the following sections,
but first we return to the general case of $k^{\text{th}}$-order
stationarity.

The $k^{\text{th}}$-order auto-cumulant function is defined to be
the function of $k-1$ variables given by
\begin{equation}
\label{eqn:cumulant} C(x_1,\ldots,x_{k-1})=
\sum_{(\nu_1,\ldots,\nu_p)}(-1)^{p-1}(p-1)!\,\mu_{\nu_1}\cdots\mu_{\nu_p}
\end{equation}
where the sum is over all partitions $(\nu_1,\ldots,\nu_p)$ of
$\{0,x_1,\ldots,x_{k-1}\}$ and $\mu_{\nu_j}=\E\left[\prod_{x_i\in
\nu_j}X_{x_i}\right]$.  For instance, for $k=3$ we have
\begin{equation}
\label{eqn:c3}
\begin{split}
C(x_1,x_2)&=\E(X_0X_{x_1}X_{x_2})-\E(X_0X_{x_1})\E(X_{x_2})
-\E(X_0X_{x_2})\E(X_{x_1})\\
&\qquad-\E(X_{x_1}X_{x_2})\E(X_0) +2\E(X_0)\E(X_{x_1})\E(X_{x_2})\\
&=\E[(X_t-\mu)(X_{t+x_1}-\mu)(X_{t+x_2}-\mu)]
\end{split}
\end{equation}
The last equality in (\ref{eqn:c3}) does not generalize for $k>3$,
that is $C$ is not the same as the more obvious generalization
\[
\widetilde{C}(x_1,\ldots,x_{k-1})=\E\left[(X_t-\mu)(X_{t+x_1}-\mu)\cdots(X_{t+x_{k-1}}-\mu)\right]
\]
For justification of the use of $C$ as opposed to $\widetilde{C}$,
the reader is referred to section 7 ``moments or cumulants?'' of
\cite{brillinger65} and section 2.3.8 ``Why Cumulant Spectra and not
Moment Spectra?'' of \cite{Nikias93b}. Nonetheless, we can see from
(\ref{eqn:cumulant}) that $C$ posses the same symmetries as
$\widetilde{C}$.  The symmetries of the general cumulant function
will be explored in subsequent sections.

The $k^{\text{th}}$-order spectral density, $f(\bm{\omega})$, is
defined as the continuous (resp. discrete) Fourier transform of the
$k^{\text{th}}$-order cumulant function, i.e.
\begin{equation*}
  f(\bm{\omega})=\frac{1}{(2\pi)^{k-1}}
  \sum_{\bm{t}\in\bZ^{k-1}}C(\bm{t})e^{-i\bm{t}\cdot\bm{\omega}}.
\end{equation*}
The second-order spectral density has many uses; see section 1.10,
``Time Series Analysis: Use of Spectral Analysis in Practice'', of
\cite{priestly81v1}.  Higher-order spectra are less frequently used,
but still many applications exist---most notably is the use of the
bispectrum in constructing linearity tests of time series
\cite{gabr80,Hinich82}. Other applications of general higher-order
spectra can be found in \cite{Nikias93b,Nikias93}.

Estimation of higher-order spectra is typically done with lag-window
estimators of the form
\begin{equation*}
\label{eq:fhat}
  \hat{f}(\bm{\omega})=\frac{1}{(2\pi)^{k-1}}\sum_{\|\bm{t}\|<N}\lambda(\bm{t}/M)\widehat{C}(\bm{t})e^{-i\bm{t}\cdot\bm{\omega}}
\end{equation*}
where $\lambda$ is a lag-window function whose properties will be
described below; $M$ is a bandwidth or smoothing parameter;
$\widehat{C}$ is an estimate of $C$ produced by estimating the
expectations $\mu_{\nu_i}$ in (\ref{eqn:cumulant}) by sample means
$\hat{\mu}_{\nu_i}$.  Indeed, estimation of $k^{\text{th}}$-order
spectra for large $k$ is plagued by the curse of dimensionality, so
accurate estimation of high-order spectra typically requires a large
sample size.

The lag-window function $\lambda$ can have varying restrictions, but
three consistent assumptions are:
\begin{problist}
  \item[(a)]  $\lambda(\bm{0})=1$;
  \item[(b)] $\|\lambda\|_{L_2}<\infty$;
  \item[(c)] $\lambda$ inherits the symmetries of $C$.
\end{problist}
It is easy to find functions to satisfy conditions (a) and (b).
However, in order to satisfy condition (c), we must first determine
the symmetries of $C$ and then construct a function that possesses
these same symmetries.  For $k=2$ the problem is trivial since the
only symmetry condition on $C(x)$ is $C(x)=C(-x)$.  We can easily
impose evenness on $\lambda$ by constructing
$\tilde{\lambda}=h\left(\lambda(x),\lambda(-x)\right)$ where $h$ is
any symmetric function (of two variables).  For instance, if
$h(x,y)=x+y$ or $h(x,y)=xy$, then $\tilde{\lambda}$ will be forced
to be even.  The solution to the general problem will be addressed
in the following sections.

\section{The Case $k=3$}

Equation \eqref{eqn:bivariatesym} in the previous section exposes
the six symmetries of $C(x,y)$, and a general technique to arrive at
these symmetries is addressed below.  We start with $S_3$, the set
of permutations of three labels.  Take for instance the permutation
$\sigma\in S_3$ that transposes the first two labels and fixes the
third. This permutation is denoted as $\sigma=(12)(3)$, or more
simply, $\sigma=(12)$.  All six permutations of $S_3$ can be written
down similarly; they constitute the set
$\{e,(12),(23),(13),(123),(132)\}$ where $e$ is the identity
permutation that doesn't move any labels.

A methodical procedure is now presented that produces each of the
symmetries of $C(x,y)$ in (\ref{eqn:bivariatesym}).  Under the
commutativity property of multiplication, the product
$X_tX_{t+x}X_{t+y}$ remains the same under any permutation of the
three variables, so for each permutations in $S_3$, we permute the
variables, according to the permutation, then adjust $t$ so that the
first variable has index $t$.  Then the symmetry condition can be
read off from the last two variables.  Say for example we take the
permutation $\sigma=(12)$, then we have
\begin{equation}
\label{eqn:example1}
\E\left[X_tX_{t+x}X_{t+y}\right]\stackrel{\sigma}{\longrightarrow}
\E\left[X_{t+x}X_tX_{t+y}\right]\longrightarrow
\E\left[X_tX_{t-x}X_{t+y-x}\right]
\end{equation}
which corresponds to the condition $C(x,y)=C(-x,y-x)$.  This process
can be simplified slightly by only writing down the indices and
introducing the function $\psi:\bR^3\rightarrow \bR^2$ given by
$\psi(a,b,c)\mapsto(b-a,c-a)$.  Then the example done in
(\ref{eqn:example1}) is simplified to
\[
(0,x,y)\stackrel{\sigma}{\longrightarrow}
(x,0,y)\stackrel{\psi}{\longrightarrow}(-x,y-x)
\]
Now suppose we take $2\times 2$ matrix resulting from the
coefficients of the $x$ and $y$ variables in each coordinate. This
induces a mapping from each permutation into the set of $2\times 2$
matrices over $\bR$.  For instance in the above example we would
have the matrix $\left(\begin{matrix} -1 & 0\\ -1 &
1\end{matrix}\right)$ corresponding to the permutation
$\sigma=(12)$.  This also establishes a correspondence between $S_3$
and identities on $C(x,y)$.  Writing down this correspondence for
each permutation in $S_3$ gives

\begin{center}
\begin{minipage}[H]{0.5\linewidth}
\begin{equation*}
\begin{split}
e&\longleftrightarrow\ \left(\begin{matrix} 1 & 0\\ 0 &
1\end{matrix}\right)\ \longleftrightarrow C(x,y)\\
(12)& \longleftrightarrow \left(\begin{matrix} -1 & 0\\ -1 &
1\end{matrix}\right)\longleftrightarrow C(-x,y-x)\\
(23)& \longleftrightarrow\ \left(\begin{matrix} 0 & 1\\ 1 &
0\end{matrix}\right)\ \longleftrightarrow C(y,x)\\
\end{split}
\end{equation*}
\end{minipage}
\begin{minipage}[H]{0.49\linewidth}
\begin{equation}
\label{eqn:matrices}
\begin{split}
(13)& \longleftrightarrow \left(\begin{matrix} 1 & -1\\ 0 &
-1\end{matrix}\right)\longleftrightarrow C(x-y,-y)\\
(123)& \longleftrightarrow \left(\begin{matrix} 0 & -1\\ 1 &
-1\end{matrix}\right)\longleftrightarrow C(-y,x-y)\\
(132)& \longleftrightarrow \left(\begin{matrix} -1 & 1\\ -1 &
0\end{matrix}\right)\longleftrightarrow C(y-x,-x)
\end{split}
\end{equation}
\end{minipage}
\end{center}
We see that each of these matrices is invertible, so in fact there
is a map, call it $\rho$, from the symmetric group $S_3$ to the
general linear group $\mathrm{GL}_2(\bR)$, the group of invertible
$2\times2$ matrices over $\bR$. Amazingly, there is compatibility in
the multiplication, namely if we take two permutations $\sigma$ and
$\tau$ and construct the composite permutation $\gamma=\sigma\tau$,
then in fact $\rho(\gamma)=\rho(\sigma)\rho(\tau)$.  In other words
$\rho$ is a group homomorphism, and since it maps a group to the
general linear group of some dimension, it is called a group
representation. For example, if we take $\sigma=(12)$ and
$\tau=(13)$, then $\sigma\tau$, going from right to left, is
equivalent to the permutation $\gamma=(132)$, and indeed
\[
\rho(\gamma)=\left(\begin{matrix} -1 & 1\\ -1 &
0\end{matrix}\right)=\left(\begin{matrix} -1 & 0\\ -1 &
1\end{matrix}\right)\left(\begin{matrix} 1 & -1\\ 0 &
-1\end{matrix}\right)=\rho(\sigma)\rho(\tau)
\]
This is more than just a group homomorphism---it is an injective
group homomorphism since each matrix in the image is distinct.  In
terms of representation theory, this is called a faithful group
representation.  This means there is no reduction of the group
structure when viewing the group in terms of its image in the matrix
group.  A simple example of a non-faithful representation is the
trivial representation that takes every element of the group to the
number 1 (which belongs to $\mathrm{GL}_1(\bR)$); in this case all
group structure is lost in the image of the representation.

Matrices similar to those in (\ref{eqn:matrices}) appear in the
first chapters of most books on group representations; for instance
one can find a similar set of the matrices on page 7 of Diaconis'
book \cite{Diaconis88}.  Yet the derivation of these matrices
typically arise in a much different context that proceeds as
follows. Let $V=\{(a,b,c)\in\bR^3| a+b+c=0\}$, a two dimensional
subspace of $\bR^3$, and consider the basis
$\{b_1=(1,-1,0),b_2=(1,0,-1)\}$. When the coordinates of $b_1$ and
$b_2$ are permuted by a permutation in $S_3$, the resulting vector
still lies in $V$ since the coordinates still sum to 0.  The
permuted basis vectors are rewritten as a linear combination of the
original basis vectors, and the coefficients are extracted to form a
$2\times 2$ matrix.  For example, $\sigma=(12)$ takes $b_1$ to
$-b_1$ and $b_2$ to $b_2-b_1$ giving the matrix correspondence for
$\sigma$ as the one in (\ref{eqn:matrices}).  In fact this
representation is identical to the one described earlier, even
though the general procedures giving rise to the representations are
completely different.

We come back to the two questions posed in the previous section by
identifying each identity in (\ref{eqn:bivariatesym}) with an
element of $S_3$ via the representation above.  If we take the
identities represented by the permutations $\sigma=(12)$ and
$\tau=(23)$, then we automatically have the identity corresponding
to $\sigma\tau=(123)$, i.e. given $C(x,y)=C(-x,y-x)$ and
$C(x,y)=C(y,x)$, then $C(x,y)=C(-y,x-y)$ follows automatically by
applying the first identity to the second.  But more identities can
be produced from just $\sigma$ and $\tau$.  For instance,
$(13)=\tau\sigma\tau$, so if we apply the identity $C(x,y)=C(y,x)$
again to $C(-y,x-y)$ we pick up $C(x-y,-y)$, i.e. the identity
corresponding to the permutation $(13)$.  Actually every identity
can be produced from just the two identities corresponding to the
permutations $\sigma$ and $\tau$. This is because these two
transpositions generate the entire group $S_3$. Therefore we can
conclude that requiring
\[
C(x,y)=C(y,x)=C(-x,y-x)
\]
is equivalent to requiring the entire string of equalities in
(\ref{eqn:bivariatesym}), and in general, any set of equations is
sufficient as long as the corresponding permutations generate all of
$S_3$. So we have established a three-way correspondence between
elements of $S_3$, matrices in $\mathrm{GL}_2(\bR)$, and symmetries
of $C(x,y)$.

There have been several attempts at constructing lag-window
functions possessing the same symmetries as in
(\ref{eqn:bivariatesym}), i.e. constructing a lag-window
$\lambda(x,y)$ that satisfies
\begin{equation}
\label{eqn:bivariatesym2}
  \lambda(x,y)=\lambda(-x,y-x)=\lambda(y,x)=\lambda(x-y,-y)=\lambda(-y,x-y)=\lambda(y-x,-x).
\end{equation}
For instance, the procedures in \cite{rao84, rao03, Nikias87,
Sasaki75,Lu96} all produce bivariate lag-windows as constructed from
univariate lag-windows. Specifically, given any even function
$\lambda(x)$,
\begin{equation}
\label{eqn:decomp} \lambda(x,y)=\lambda(x)\lambda(y)\lambda(x-y)
\end{equation}
will satisfy the requirements of (\ref{eqn:bivariatesym2}).  However
this class of functions is too restrictive; we will see below that
the popular Gabr-Rao optimal lag-window is excluded from this
construction.   A much different approach to constructing bivariate
lag-windows is considered in \cite{Yao01}, but their focus is not on
symmetrizing the lag-window function.

 The optimal lag-window, $\lambda_{\text{opt}}$, introduced in
\cite{rao84}, is optimal in the sense that the variance of the
bispectral estimators using this lag-window is least among a certain
class of lag-windows; refer to Theorem 3 in the next section for
generalization of this property.  $\lambda_{\text{opt}}$ is defined
by the inverse Fourier transform of the optimal kernel
$\Lambda_{\text{opt}}$ given by
\begin{equation}
\label{eqn:GR}
  \Lambda_{\text{opt}}(\omega_1,\omega_2)=
\begin{cases}
\frac{\sqrt{3}}{\pi^3}\left(1-\frac{1}{\pi^2}\left(\omega_1^2+\omega_2^2+\omega_1\omega_2\right)\right),
& \text{if } \omega_1^2+\omega_2^2+\omega_1\omega_2\le\pi^2\\
  0, &\text{otherwise}
\end{cases}
\end{equation}
This kernel, as with every kernel, is unique up to scale, i.e. the
kernel $\Lambda_{\text{opt}}(x,y)$ is equivalent to the kernel
$h\Lambda_{\text{opt}}(hx,hy)$ for any $h>0$.  Only an approximation
of the inverse Fourier transform of $\Lambda_{\text{opt}}$ is
provided in \cite{rao84}, but in \cite{Saito85},
$\lambda_{\text{opt}}$ is shown to be
\[
\lambda_{\text{opt}}(\tau_1,\tau_2)=\frac{8}{\alpha(\tau_1,\tau_2)^2}
J_2(\alpha(\tau_1,\tau_2))
\]
where
\[
\alpha(x,y)=\frac{2\pi}{\sqrt{3}}\sqrt{x^2-xy+y^2}
\]
and $J_2$ is the second-order Bessel function of the first kind,
i.e.
\[
J_2(x)=\sum_{\ell=0}^\infty
\frac{(-1)^\ell}{2^{2\ell+2}\ell!(2+\ell)!}x^{2\ell+2}.
\]
We now ask the question: does $\lambda_{\text{opt}}(x,y)$ admit a
decomposition of the form (\ref{eqn:decomp})?  It is immediately
true that no such decomposition exists for $\alpha(x,y)$ since
otherwise $\alpha(0,0)=0$ would imply $\lambda(0)=0$, but then
$\alpha(1,1)$ would be forced to be zero (when it's not).  We extend
the answer to $\lambda_{\text{opt}}(x,y)$ with the following
theorem.

\begin{theorem}
Suppose $\lambda(x,y)$ is a continuous function that vanishes on a
set of measure zero, and suppose $\lambda$ does vanish on some
smooth nonlinear curve $\mathcal{C}\subset\bR^2$.   Additionally
suppose $\lambda(0,0)\not=0$.  Then there is no function $f$ such
that $\lambda(x,y)=f(x)f(y)f(x-y)$.
\end{theorem}
\begin{proof}
Suppose there is a function $f$ such that
$\lambda(x,y)=f(x)f(y)f(x-y)$.  If $f(x)=0$ on some interval $I$,
then $\lambda(x,y)=0$ for every $(x,y)\in I\times I$ contradicting
the assumption that $\lambda$ vanishes on a set of measure zero.  We
will now show that $f(x)$ is zero on some interval.  Since
$\lambda(0,0)=f(0)^3\not=0$, we get that $f(0)\not=0$.  Also,
$\lambda(x,x)=f(x)^2f(0)$ is continuous so $f(x)^2$ is continuous.
Let $(x^*,y^*)$ be any point on $C$, then
$f(x^*)f(y^*)f(x^*-y^*)=0$, so at least one of $f(x^*)$, $f(y^*)$,
$f(x^*-y^*)$ must be zero.  As $(x^*,y^*)$ moves continuously along
$\mathcal{C}$, the smooth and nonlinear properties of $\mathcal{C}$
guarantee that $x^*$, $y^*$, $x^*-y^*$ each simultaneously sweep out
some interval of points.  Continuity of $f^2$ implies the continuity
of $|f(x^*)|$, $|f(y^*)|$, and $|f(x^*-y^*)|$, but since their
product is always zero on $\mathcal{C}$, there must be some interval
in which $f$ is zero.

\end{proof}

Since $\lambda_{\text{opt}}$ is continuous, vanishes on a set of
measure zero, and vanishes on the ellipse $x^2-xy+y^2=2$, we have
the following corollary.
\begin{corollary}
The optimal lag-window, $\lambda_{\text{opt}}$, cannot be written as
$\lambda_{\text{opt}}(x,y)=f(x)f(y)f(x-y)$ for any choice of
function $f$.
\end{corollary}

Now we will construct a lag-window satisfying
(\ref{eqn:bivariatesym2}) from any bivariate function $f(x,y)$.
We have already seen a symmetrizing technique to create an even
function in the previous section, and in taking a similar approach
to this, we symmetrize $f(x,y)$ by
\begin{equation}
\label{eqn:fsym} \tilde{f}(x,y)=h(f(x,y), f(-x,y-x),f(y,x),
f(x-y,-y),f(-y,x-y),f(y-x,-x))
\end{equation}
where $h$ is any symmetric function of its six arguments; for
instance $h$ could be a power mean like the arithmetic mean or
geometric mean.  The above construction can be made slightly more
general by replacing $f(x,y)$ in (\ref{eqn:fsym}) with $k(f(x,y))$
where $k$ is any (univariate) function.  On the other hand if we
take just the first coordinates in each $f$ in the above formula for
$\tilde{f}$, this also produces a function that satisfies the
required symmetries. That is, if we let
\[
\tilde{g}(x,y)= h(g(x),g(-x),g(y),g(x-y),g(-y),g(y-x))
\]
then $\tilde{g}(x,y)$ is another symmetrization built up from a
generic univariate function $g$, and if $g$ is assumed to be an even
function, then $\tilde{g}(x,y)$ reduces to
\[
\tilde{\tilde{g}}(x,y)=\tilde{h}\left(g(x),g(y),g(x-y)\right)
\]
where $\tilde{h}$ is any symmetric function of three variables.
Therefore we automatically get the construction (\ref{eqn:decomp})
with the special case $\tilde{h}(x,y,z)=xyz$.

To prove that these symmetrizations $\tilde{f}$ and $\tilde{g}$ do
satisfy the required symmetries, we return to the connection with
group representations and also introduce the concept of group
actions.  If we fix any permutation in $S_3$ and multiply it with
all the elements of $S_3$, the resulting action just permutes the
six elements of $S_3$; this is a type of group action where $S_3$
``acts on itself''. To see that it is just a permutation of the
group elements, suppose $\sigma\tau=\sigma\gamma$, then multiplying
both sides by $\sigma^{-1}$ gives $\tau=\gamma$.  Suppose we wish to
test a symmetry like $\tilde{f}(x,y)=\tilde{f}(y,x)$. We see from
(\ref{eqn:matrices}) that this condition corresponds to the
permutation $(23)$, and multiplying $(23)$ with each element of
$S_3$ just reorders the elements of $S_3$, i.e.
\begin{equation}
\label{eqn:action}
(23)\cdot\{e,(12),(23),(13),(123),(132)\}=\{(23),(132),e,(123),(13),(12)\}
\end{equation}
The righthand side corresponds to the equation
\[
\tilde{f}(y,x)=h(f(y,x), f(y-x,-x),f(x,y),
f(-y,x-y),f(x-y,-y),f(-x,y-x))
\]
and since $h$ is a symmetric function, this is equivalent to the
original equation $\tilde{f}(x,y)$.  In general, since each identity
required just corresponds to a permutations of the arguments of $h$,
each identity will be satisfied. To justify the symmetries of
$\tilde{\tilde{f}}$, we note that taking the first coordinate
corresponds to multiplying the matrix representation with the vector
$(1\ \ 0)$.  So instead of identifying each permutation with the
matrix in (\ref{eqn:matrices}), it is identified by the $1\times 2$
vector that is the top row of each matrix.  Now if we wish to show
$\tilde{g}(x,y)=\tilde{g}(y,x)$, we see that this is equivalent to
multiplying each $2\times 1$ vector by the matrix
$\left(\begin{matrix} 0 & 1\\ 1 & 0\end{matrix}\right)$.  So again
we see that this just induces a permutation of the elements, i.e.
\[
\left(\begin{matrix} 1 & 0\end{matrix}\right)\cdot \left\{
\left(\begin{matrix} 1 & 0\\ 0 & 1\end{matrix}\right),
\left(\begin{matrix} -1 & 0\\ -1 & 1\end{matrix}\right),
\left(\begin{matrix} 0 & 1\\ 1 & 0\end{matrix}\right),
\left(\begin{matrix} 1 & -1\\ 0 & -1\end{matrix}\right),
\left(\begin{matrix} 0 & -1\\ 1 & -1\end{matrix}\right),
\left(\begin{matrix} -1 & 1\\ -1 & 0\end{matrix}\right) \right\}
\cdot \left(\begin{matrix} 0 & 1\\ 1 & 0\end{matrix}\right)
\]
is just a permutation of the elements
\[
\left(\begin{matrix} 1 & 0\end{matrix}\right)\cdot \left\{
\left(\begin{matrix} 1 & 0\\ 0 & 1\end{matrix}\right),
\left(\begin{matrix} -1 & 0\\ -1 & 1\end{matrix}\right),
\left(\begin{matrix} 0 & 1\\ 1 & 0\end{matrix}\right),
\left(\begin{matrix} 1 & -1\\ 0 & -1\end{matrix}\right),
\left(\begin{matrix} 0 & -1\\ 1 & -1\end{matrix}\right),
\left(\begin{matrix} -1 & 1\\ -1 & 0\end{matrix}\right) \right\}
\]
 So we see the same theory as
before applies as every identity is associated with a permutation
the arguments of $\tilde{h}$ leaving the value fixed.  There is
nothing special about the vector $(1\ \ 0)$, and different vectors
(or more generally $n\times 2$ matrices) lead to similar formulas.
For instance, using the vector $(1\ \ 1)$, we see that the function
\[
g(x+y)g(x-2y)g(y-2x)
\]
also satisfies the required symmetries.

\section{The General Case}

Here we move from the symmetric group of degree three to the
symmetric group of general degree $k$.  It is well known that any
permutation can be written as a product of cycles and every cycle is
a product of transpositions (2-cycles), therefore the permutations
can be generated by just the transpositions.  But more is true:
every permutation is the product of transpositions of the form $(i
,i+1)$ for $i=1,\ldots,k-1$.  To see this, suppose $i<j$, then
\[
(i,j+1)=(j,j+1)(i,j)(j,j+1)
\]
So any transposition, after applying the above formula enough times,
can be reduced to a product of transpositions of consecutive labels.
Thus the group $S_k$ of $k!$ permutations can be represented by a
much smaller subset of only $k-1$ transpositions.

Now a representation of $S_k$ is constructed for general $k$.  We
generalize the $\psi$ function to
$\psi:\bR^k\longrightarrow\bR^{k-1}$ given by
\[
(x_1,x_2,\ldots,x_k)\mapsto(x_2-x_1,\ldots,x_k-x_1)
\]
and the representation on $S_k$ is produced by composing the
permutation with the $\psi$ function acted on
$(0,x_1,\ldots,x_{k-1})$, then extracting the $(k-1)$-dimensional
square matrix of coefficients from the image in the natural way.  As
an example, let $\sigma=(1234)\in S_4$, then
\[
(0,x,y,z)\stackrel{\sigma}{\longrightarrow}
(z,0,x,y)\stackrel{\psi}{\longrightarrow}(-z,x-z,y-z)
\]
So we have the correspondence
\[
(1234)\longleftrightarrow\left(\begin{matrix} 0 & 0 & -1\\ 1 & 0 &
-1\\0 & 1 & -1\end{matrix}\right)
\]
This process clearly induces a map from $S_k$ to $(k-1)$-dimensional
square matrices, but it is not clear a priori that the matrices are
invertible or if this operation establishes a group homomorphism.
The following theorem proves this operation is indeed a group
representation, and moreover, is a faithful group representation.
\begin{theorem}
  The map, referred to now as $\rho$, from $S_k$ to $(k-1)$-dimensional square matrices
  described above is a faithful group representation of dimension
  $k-1$.
\end{theorem}
\begin{proof}
Let $P_\sigma$ be the $k\times k$ permutation matrix associated to
$\sigma$, i.e. the $(i,j)$-entry is 1 if $\sigma$ maps $i$ to $j$
and 0 otherwise.  It is well known that permutation matrices are
invertible ($P_{\sigma}^{-1}=P_{\sigma^{-1}}$), in fact the map
$\sigma\mapsto P_\sigma$ is an $k$-dimensional group representation
of $S_k$. Let $A$ be the matrix with ones in the first column and
zeros elsewhere. Then the above mapping can be described by
associated the permutation $\sigma$ with $P_\sigma$ then the
operation of the $\psi$ function is represented by subtracting each
row of $P_\sigma$ by the first row of $P_\sigma$ and then discarding
the top row and leftmost column. Hence the representation $\rho$ is
described compactly by the map
\begin{equation}
\label{eqn:compact} \sigma\mapsto \left\langle
P_\sigma-AP_\sigma\right\rangle =\left\langle
(I-A)P_\sigma\right\rangle
\end{equation}
where the notation $\langle B\rangle$ denotes truncation of the top
row and leftmost column of the matrix $B$.  From basic block
multiplication of matrices, we see that
\begin{equation}
\label{eqn:Bs} \langle B_1\rangle\langle B_2\rangle=\langle
B_1B_2\rangle
\end{equation}
 for any
two matrices $B_1$ and $B_2$ the multiply compatibly.  Also from
block multiplication,  given any matrix $B$,
\begin{equation}
\label{eqn:for2} \left\langle B (I-A)\right\rangle=\left\langle
B\right\rangle
\end{equation}
as long as the multiplication is compatible.  We first show that the
image of a permutation under $\rho$ is invertible, i.e. an element
of $\mathrm{GL}_{k-1}(\bR)$, by exhibiting an inverse.
Specifically, we show that the inverse of the matrix $\left\langle
(I-A)P_\sigma\right\rangle$ is $\left\langle
(I-A)P_\sigma^{-1}\right\rangle$.  Applying properties
(\ref{eqn:Bs}) and (\ref{eqn:for2}) shows
\[
\begin{split}
\left\langle (I-A)P_\sigma\right\rangle \left\langle
(I-A)P_\sigma^{-1}\right\rangle &= \left\langle I-A\right\rangle
\left\langle P_\sigma (I-A)\right\rangle
\left\langle P_\sigma^{-1}\right\rangle\\
&=\left\langle I-A\right\rangle \left\langle P_\sigma\right\rangle
\left\langle P_\sigma^{-1}\right\rangle\\
&=\left\langle I-A\right\rangle\\
&=\left\langle I\right\rangle.
\end{split}
\]

In order to prove $\rho$ is a representation, it must be shown that
for any $\sigma,\tau\in S_k$,
$\rho(\sigma)\rho(\tau)=\rho(\sigma\tau)$. This is equivalent to
\begin{equation}
\label{eqn:eq1} \left\langle (I-A)P_\sigma\right\rangle \left\langle
(I-A)P_\tau\right\rangle= \left\langle (I-A)P_\sigma
P_\tau\right\rangle
\end{equation}
By multiplying both sides on the right by $\langle
P_\tau^{-1}\rangle$ and making free use of (\ref{eqn:Bs}), equation
(\ref{eqn:eq1}) becomes
\begin{equation}
\label{eqn:id1} \left\langle (I-A)P_\sigma (I-A)\right\rangle=
\left\langle (I-A)P_\sigma\right\rangle
\end{equation}
(\ref{eqn:id1}) is clearly true by noting  $\left\langle P_\sigma
(I-A)\right\rangle= \left\langle P_\sigma\right\rangle$.  Therefore
$\rho$ is indeed a group representation.

The final task is to show $\rho$ is faithful, i.e. $\rho$ is
injective. One easy way to prove this is to show that the only
permutation being mapped to $\langle I\rangle$ is $e$, the identity
permutation. In the language of group theory, this says that we need
to show the kernel of $\rho$ is the identity.  So suppose
$\left\langle (I-A)P_\sigma\right\rangle = \langle I\rangle$, then
by multiplying both sides by $\langle P_\sigma^{-1}\rangle$ gives
$\langle I\rangle =\langle P_\sigma^{-1}\rangle=\langle
P_{\sigma^{-1}}\rangle$.  This says the permutation $\sigma^{-1}$
fixes the labels 2 through $k$, so it must also fix label 1.  Hence
$\sigma^{-1}=e=\sigma$, and the only permutation in the kernel of
$\rho$ is $e$.

\end{proof}

So now that we have established $\rho$ is a true group
representation, we can extend the results from the previous section
to general $k$.  For instance if $\rho(\sigma)=(a_{i,j})$ for some
$\sigma\in S_k$, then the acf of $k-1$ variables satisfies the
identity
\begin{equation}
\label{eqn:sym} C(x_1,\ldots,x_{k-1})=C\left(\sum_{j=1}^{k-1}
a_{1,j}x_j\, ,\, \ldots\, ,\sum_{j=1}^{k-1} a_{k-1,j}x_j \right),
\end{equation}
Also, since the transpositions $(12), (23), \ldots, (k-1,k)$
generate $S_k$, all of the symmetries of $C$ can be described by
(\ref{eqn:sym}) with just these permutations.  If we consider the
permutations that fix the first label, then we see that the
identities in (\ref{eqn:sym}) induced by these $(k-1)!$ permutations
are just permutations of the arguments.  Therefore $C$, in
particular, is a symmetric function.

Now with the results of Theorem 2, we can construct lag-windows
functions that satisfy all of the required symmetries.
\begin{corollary}
Given any function $f:\bR^{k-1}\longrightarrow \bR$, a
symmetrization of $f$ is constructed by starting with any symmetric
function, $h$, of $n!$ variables and plugging-in $f_\sigma$ into the
variables of $h$ (using every $\sigma\in S_k$) where, like in
(\ref{eqn:sym}),
\[
f_\sigma(x_1,\ldots,x_{k-1})=f\left(\sum_{j=1}^{k-1} a_{1,j}x_j\,
,\, \ldots\, ,\sum_{j=1}^{k-1} a_{k-1,j}x_j \right).
\]
\end{corollary}
\section{The Fourier Transform}

We now come back to the motivation of this theory which is to
construct higher-order spectral density estimates using lag-windows
as in (\ref{eq:fhat}), i.e. estimates of the form
\[
  \hat{f}(\bm{\omega})=\frac{1}{(2\pi)^{k-1}}\sum_{\|\bm{t}\|<N}\lambda(\bm{t}/M)\widehat{C}(\bm{t})e^{-i\bm{t}\cdot\bm{\omega}}
\]
There is an equivalent expression to this estimator in the frequency
domain given by
\begin{equation}
\label{eqn:fhat2} \hat{f}(\bm{\omega})=\Lambda\ast
P(\bm{\omega})=\int_{\bR^{k-1}}\Lambda(\bm{t}/M)
P(\bm{\omega}-\bm{t})\,d\bm{t}
\end{equation}
where $\Lambda$ is the Fourier transform of $\lambda$ and $P$ is the
$(k-1)^{\text{th}}$ order periodogram; namely,
\[
\Lambda(\bm{t})=\int_{\bR^{k-1}}
\lambda(\bm{t})e^{-i\bm{\omega}\cdot \bm{t}}\,d\bm{t}
\]
and
\[
P(\bm{\omega})=\frac{1}{(2\pi)^{k-1}}\sum_{ \bm{t}\in
\bZ^{k-1}}\widehat{C}(\bm{t})e^{-i\bm{t}\cdot\bm{\omega}}
\]
If we were to use the less practical estimator given in
(\ref{eqn:fhat2}), then we would need a kernel $\Lambda$ instead of
a lag-window $\lambda$, so an understanding of the symmetries of
$\Lambda$ is also important. To begin, we start by considering an
example in the case $k=3$.  The permutation $\sigma=(12)$
corresponds to the identity $\lambda(x,y)=\lambda(-x,y-x)$, so we
consider what the substitutions $\omega_1\mapsto-\omega_1$ and
$\omega_2\mapsto\omega_2-\omega_1$ do to $\Lambda$;
\[
\begin{split}
\Lambda(-\omega_1,\omega_2-\omega_1) &= \int_{-\infty}^\infty
\int_{-\infty}^\infty \lambda(x,y)
e^{-ix(-\omega_1)}e^{-iy(\omega_2-\omega_1)}\,dx\,dy\\
&=\int_{-\infty}^\infty \int_{-\infty}^\infty \lambda(x,y)
e^{-i\omega_1(-x-y)}e^{-i\omega_2(y)}\,dx\,dy.
\end{split}
\]
After the simplification in the second line, the actions of the
Fourier transform begin to unfold, i.e. the exponential kernel of
the Fourier transform induced the transposed representation
$x\mapsto -x-y$ and $y\mapsto y$.  Now if we start with this
transposed identity, i.e. $\omega_1\mapsto-\omega_1-\omega_2$ and
$\omega_2\mapsto\omega_2$, we find this produces a symmetry of
$\Lambda$:
\[
\begin{split}
\Lambda(-\omega_1-\omega_2,\omega_2) &= \int_{-\infty}^\infty
\int_{-\infty}^\infty \lambda(x,y)
e^{-ix(-\omega_1-\omega_2)}e^{-iy\omega_2}\,dx\,dy\\
&=\int_{-\infty}^\infty \int_{-\infty}^\infty \lambda(x,y)
e^{-i\omega_1(-x)}e^{-i\omega_2(y-x)}\,dx\,dy\\
&=\int_{-\infty}^\infty \int_{-\infty}^\infty \lambda(-x,y-x)
e^{-i\omega_1x}e^{-i\omega_2y}\,dx\,dy\\
&=\Lambda(\omega_1,\omega_2)
\end{split}
\]

Thus we have a pseudo-representation $\tilde{\rho}$ which maps
$\sigma\in S_k$ to $\rho(\sigma)'$, the transpose of $\rho(\sigma)$.
For $\sigma, \tau\in S_k$, This representation satisfies
\[
\tilde{\rho}(\sigma\tau)=\rho(\sigma\tau)'=\left(\rho(\sigma)\rho(\tau)\right)'=\rho(\tau)'\rho(\sigma)'=\tilde{\rho}(\tau)\tilde{\rho}(\sigma)
\]
which is close to the requirement of a group representation, just
reversed.  However if we consider multiplication of two permutations
in $S_k$ in the opposite order, from left to right, then in fact
$\tilde{\rho}$ is a true group representation with the same
properties as $\rho$.  For example the multiplication $(12)(13)$
from left to right becomes $(123)$ and not $(132)$ as in the right
to left situation.  Therefore we see that $\Lambda$ possesses the
``transposed symmetries'' as those of $\lambda$, i.e. if
$\tilde{\rho}(\sigma)=(b_{i,j})$ for some $\sigma\in S_k$, then
$\Lambda$ satisfies
\begin{equation}
\label{eqn:sym2}
\Lambda(\omega_1,\ldots,\omega_{k-1})=\Lambda\left(\sum_{j=1}^{n-1}
b_{1,j}\omega_j\, ,\, \ldots\, ,\sum_{j=1}^{n-1} b_{n-1,j}\omega_j
\right).
\end{equation}
In particular, like in the case of the lag-window $\lambda$,
$\Lambda$ is a symmetric function.

As an example, the symmetries for $k=3$ are given by
\[
\Lambda(x,y)=\Lambda(-x-y,y)=\Lambda(y,x)=\Lambda(x,-x-y)=\Lambda(y,-x-y)=\Lambda(-x-y,y)
\]
A construction of the bivariate kernels from univariate kernels
similar to \eqref{eqn:decomp} is easily derived using the theory
described to become $\Lambda(x,y)=\Lambda(x)\Lambda(y)\Lambda(-x-y)$
where the only requirement on the univariate $\Lambda(x)$ function
is nonnegativity.

 We shall now apply this theory to generalize the
Gabr-Rao optimal window, $\Lambda_{\text{opt}}$, for estimation of
the trispectrum and general $k^{\text{th}}$-order spectra. The
window must integrate to one and satisfy (\ref{eqn:sym2}) for every
$\sigma\in S_k$. Observing that the usual Gabr-Rao window is a
constant plus a homogeneous polynomial of degree 2 (inside its
support), we start with the construction
\begin{equation}
\label{eqn:exp1}
\tilde{\Lambda}_{\text{opt}}(\omega_1,\ldots,\omega_{k-1})=c+\sum_{i\leq
j} c_{ij}\,\omega_i\,\omega_j
\end{equation}
for some constants $c$ and $c_{ij}$.  Since $\Lambda_{\text{opt}}$
is a symmetric function, equation (\ref{eqn:exp1}) can be greatly
simplified to
\begin{equation}
\label{eqn:exp2}
\tilde{\Lambda}_{\text{opt}}(\omega_1,\ldots,\omega_{k-1})=c+c_1\sum_{i=1}^{k-1}
\omega_i^2 +c_2\sum_{i<j} \omega_i\,\omega_j
\end{equation}
Here we see the convenience in associating the symmetries of
$\Lambda$ to $S_k$.  We showed earlier that it is not necessary to
test all of the possible identities of $\Lambda$, just the ones that
correspond to the permutations that generate $S_k$.  Since equation
(\ref{eqn:exp2}) accounts for the symmetries corresponding to the
permutations that fix 1, all that is left is to check the identity
corresponding to the permutation $(12)$, i.e. if
$\tilde{\Lambda}_{\text{opt}}$ satisfies
\begin{equation}
\label{eqn:12}
\tilde{\Lambda}_{\text{opt}}(\omega_1,\ldots,\omega_{k-1})=\tilde{\Lambda}_{\text{opt}}
\left(-\sum_{i=1}^{k-1} \omega_i\,
,\,\omega_2\,,\,\ldots\,,\,\omega_{k-1}\right)
\end{equation}
then $\Lambda_{\text{opt}}$ satisfies all of the necessary
symmetries.  Applying (\ref{eqn:12}) to (\ref{eqn:exp2}) gives
\begin{subequations}
\begin{align}
\tilde{\Lambda}_{\text{opt}}&(\omega_1,\ldots,\omega_{k-1})=c+c_1\sum_{i=1}^{k-1}\omega_i^2+c_2\sum_{i<j}\omega_i\omega_j\\
&=c+c_1\omega_1^2 + c_2\omega_1\sum_{j=2}^{k-1}\omega_j + c_1\sum_{i=1}^{k-1}\omega_i^2+c_2\sum_{1<i<j}\omega_i\omega_j\label{second}\\
&=c+c_1\left(\sum_{i=1}^{k-1}\omega_i\right)^2 -
c_2\left(\sum_{i=1}^{k-1}\omega_i\right)\sum_{j=2}^{k-1}\omega_j +
c_1\sum_{i=1}^{k-1}\omega_i^2+c_2\sum_{1<i<j}\omega_i\omega_j\label{third}
\end{align}
\end{subequations}
Comparing \eqref{second} with \eqref{third} gives the identity
\begin{equation}
\label{eqn:identity} c_1\omega_1^2 +
c_2\omega_1\sum_{j=2}^{k-1}\omega_j=c_1\left(\sum_{i=1}^{k-1}\omega_i\right)^2
- c_2\left(\sum_{i=1}^{k-1}\omega_i\right)\sum_{j=2}^{k-1}\omega_j
\end{equation}
There are many different ways to massage \eqref{eqn:identity} into
conditions on $c_1$ and $c_2$.  Here we differentiate both sides of
\eqref{eqn:identity} with respect to $\omega_1$ producing
\[
2c_1\omega_1+c_2\sum_{i=2}^{k-1}\omega_i=2c_1\sum_{i=1}^n\omega_i
-c_2\sum_{i=2}^{k-1}\omega_i
\]
which is equivalent to
\[
2c_1\sum_{i=2}^{k-1}\omega_i=2c_2\sum_{i=2}^{k-1}\omega_i
\]
Therefore we must have $c_1=c_2$, and it is easily seen that
\eqref{eqn:identity} is satisfied under this condition.  Thus
\begin{equation}
\label{eqn:poly}
\tilde{\Lambda}_{\text{opt}}(\omega_1,\ldots,\omega_{k-1})=
\alpha\left(1-\beta\left(\sum_{i=1}^{k-1}\omega_i^2+\sum_{i<j}\omega_i\omega_j\right)\right)
\end{equation}
satisfies all of the necessary symmetries.  We now define
$\Lambda_{\text{opt}}(\bm{\omega})$ as
\begin{equation}
\label{eqn:Lam_opt}
\Lambda_{\text{opt}}(\bm{\omega})=\tilde{\Lambda}(\bm{\omega})^{+}=
\begin{cases}
  \tilde{\Lambda}(\bm{\omega}), & \text{if } \tilde{\Lambda}(\bm{\omega})\ge0\\
  0, &\text{otherwise}
\end{cases}
\end{equation}
where $\beta$ is any positive constant (kernels are unique up to
scale) and $\alpha$ is chosen such that $\Lambda_{\text{opt}}$
integrates to one.
\begin{theorem}
Let $\Lambda(\bm{\omega})$ be any nonnegative kernel that integrates
to one and satisfies all the necessary symmetries, i.e. satisfies
\eqref{eqn:sym2} for all $\sigma\in S_{k}$.  Also assume
\begin{equation}
\label{eqn:assumption}
  \int_{\bR^{k-1}}\omega_j^2\Lambda(\bm{\omega})\,d\bm{\omega}=
  \int_{\bR^{k-1}}\omega_j^2\Lambda_{\mathrm{opt}}(\bm{\omega})\,d\bm{\omega}
\end{equation}
for $j=1,\ldots,n-1$.  Then
$\|\Lambda\|_{L_2}\ge\|\Lambda_{\mathrm{opt}}\|_{L_2}$, i.e.
\[
  \int_{\bR^{k-1}}\Lambda(\bm{\omega})^2\,d\bm{\omega}\ge
  \int_{\bR^{k-1}}\Lambda_{\mathrm{opt}}(\bm{\omega})^2\,d\bm{\omega}.
\]
\end{theorem}
\begin{proof}
Let
$\Lambda(\bm{\omega})=\Lambda_{\text{opt}}(\bm{\omega})+\varepsilon(\bm{\omega})$,
then plugging this substitution into \eqref{eqn:assumption} gives
\[
\int_{\bR^{k-1}}\omega_j^2\Lambda(\bm{\omega})\,d\bm{\omega}=
\int_{\bR^{k-1}}\omega_j^2\left(\Lambda_{\text{opt}}(\bm{\omega})+\varepsilon(\bm{\omega})\right)\,d\bm{\omega}=
\int_{\bR^{k-1}}\omega_j^2\Lambda_{\text{opt}}(\bm{\omega})\,d\bm{\omega}
\]
Therefore
\begin{equation}
\label{eqn:zero1}
\int_{\bR^{k-1}}\omega_j^2\,\varepsilon(\bm{\omega})\,d\bm{\omega}=0
\end{equation}
for all $j$.  Again making use of the identity in \eqref{eqn:12}
gives
\[
\begin{split}
0&=\int_{\bR^{k-1}}\omega_1^2\,\varepsilon(\bm{\omega})\,d\bm{\omega}=
\int_{\bR^{k-1}}\bigg(\sum_{j=1}^{k-1}
\omega_j\bigg)^2\,\varepsilon(\bm{\omega})\,d\bm{\omega}\\
&=\sum_{j=1}^{k-1} \int_{\bR^{k-1}}
\omega_j^2\,\varepsilon(\bm{\omega})\,d\bm{\omega}+
2\sum_{i<j}\int_{\bR^{k-1}}\omega_i\,\omega_j\,\varepsilon(\bm{\omega})\,d\bm{\omega}\\
&=2\sum_{i<j}\int_{\bR^{k-1}}\omega_i\,\omega_j\,\varepsilon(\bm{\omega})\,d\bm{\omega}
\end{split}
\]
Therefore
\begin{equation}
\label{eqn:zero2}
\sum_{i<j}\int_{\bR^{k-1}}\omega_i\,\omega_j\,\varepsilon(\bm{\omega})\,d\bm{\omega}=0.
\end{equation}
Since $\Lambda$ and $\Lambda_{\text{opt}}$ both integrate to one,
$\varepsilon$ must integrate to zero, i.e.
\begin{equation}
\label{eqn:zero3}
\int_{\bR^{k-1}}\varepsilon(\bm{\omega})\,d\bm{\omega}=0.
\end{equation}
Computing the $L_2$-norm of $\Lambda(\bf{\omega})$ gives
\[
\int_{\bR^{k-1}}\Lambda(\bm{\omega})^2\,d\bm{\omega}=
\int_{\bR^{k-1}}\Lambda_{\text{opt}}(\bm{\omega})^2\,d\bm{\omega}+
\int_{\bR^{k-1}}\varepsilon(\bm{\omega})^2\,d\bm{\omega}+
2\int_{\bR^{k-1}}\Lambda_{\text{opt}}(\bm{\omega})\varepsilon(\bm{\omega})\,d\bm{\omega}
\]
So if we can show that last summand is nonnegative, we will have
proved the theorem.  Referring to the definition of
$\tilde{\Lambda}(\bm{\omega})$ in \eqref{eqn:poly}, we have
\begin{equation}
\label{eqn:integral}
\int_{\bR^{k-1}}\Lambda_{\text{opt}}(\bm{\omega})\varepsilon(\bm{\omega})\,d\bm{\omega}=
\int_{\bR^{k-1}}\tilde{\Lambda}_{\text{opt}}(\bm{\omega})\varepsilon(\bm{\omega})\,d\bm{\omega}-
\int\limits_{\bm{\omega}:\tilde{\Lambda}_{\text{opt}}(\bm{\omega})<0}
\tilde{\Lambda}_{\text{opt}}(\bm{\omega})\varepsilon(\bm{\omega})\,d\bm{\omega}
\end{equation}
The first integral is zero by \eqref{eqn:zero1}, \eqref{eqn:zero2},
and \eqref{eqn:zero3}.  When $\tilde{\Lambda}(\bm{\omega})<0$,
$\Lambda_{\text{opt}}(\bm{\omega})=0$ and $\Lambda(\bm{\omega})\ge
0$, so we must have $\varepsilon(\bm{\omega})\ge0$. Therefore the
second integral is less than zero making the left hand side of
\eqref{eqn:integral} nonnegative.
\end{proof}

For $k=2$, if we let $\beta=1/5$, then $\alpha=\frac{3}{4\sqrt{5}}$
and \eqref{eqn:Lam_opt} is the familiar Epanechnikov kernel
\cite{epanechnikov69}. This kernel is equivalent to the
Bartlett-Priestly kernel in \cite{priestly81v1} and corresponds to
the quadratic spectral lag-window in \cite{Andrews91}.  For $k=3$,
$\beta=1/\pi^2$, and $\alpha=\sqrt{3}/\pi^3$, \eqref{eqn:Lam_opt} is
the Gabr-Rao optimum bispectral kernel \cite{rao84}.

The set of points that satisfies
$\tilde{\Lambda}_{\text{opt}}(\bm{\omega})>0$ is equivalent to
\[
\sum_{i=1}^{k-1}\omega_i^2+\sum_{i<j}\omega_i\omega_j \le
\frac{1}{\beta}
\]
where the left hand side is a positive quadratic form since
\begin{equation}
\label{eqn:quad}
\sum_{i=1}^{k-1}\omega_i^2+\sum_{i<j}\omega_i\omega_j=\frac{1}{2}\sum_{i=1}^{k-1}\omega_i^2+
\frac{1}{2}\left(\sum_{i=1}^{k-1}\omega_i\right)^{2}\ge 0
\end{equation}
and equal to zero only if $\omega_j=0$ for all $j$.  Therefore there
always exists some transformation of coordinates to transform the
ellipsoid defined by $\tilde{\Lambda}_{\text{opt}}(\bm{\omega})>0$
to the unit sphere, i.e. a diagonalization of the quadratic form
\eqref{eqn:quad} to standard form with rank and signature $k-1$.
This transformation is easily computed for any given $k$
\cite{Beightler66}, and the Jacobian of the transformation can be
used to determine $\alpha$ for a given $\beta$ and $k$.  For
instance for $k=4$, the substitutions
\[
\begin{split}
&\omega_1\longrightarrow
u_1-\frac{u_2}{\sqrt{3}}-\frac{u_3}{\sqrt{6}}\\
&\omega_2\longrightarrow
-u_1-\frac{u_2}{\sqrt{3}}-\frac{u_3}{\sqrt{6}}\\
&\omega_3\longrightarrow \sqrt{\frac{3}{2}}u_3
\end{split}
\]
transforms
$\omega_1^2+\omega_2^2+\omega_3^2+\omega_1\omega_2+\omega_2\omega_3+\omega_1\omega_3$
into $u_1^2+u_2^2+u_3^2$.  The Jacobian of this transformation is
$\sqrt{2}$, and so $\alpha$ is computed to be
\[
\alpha=\frac{15\beta^{3/2}}{8\sqrt{2}\pi}.
\]
Polar coordinates on $\bR^{k-1}$ can be used to exactly determine
$\alpha$ for larger $k$ \cite{Folland99}.

\section{Bivariate Examples}

The bivariate optimal kernel with corresponding lag-window are
plotted below.

\begin{figure}[H]
\centerline{\includegraphics[width=2.5in]{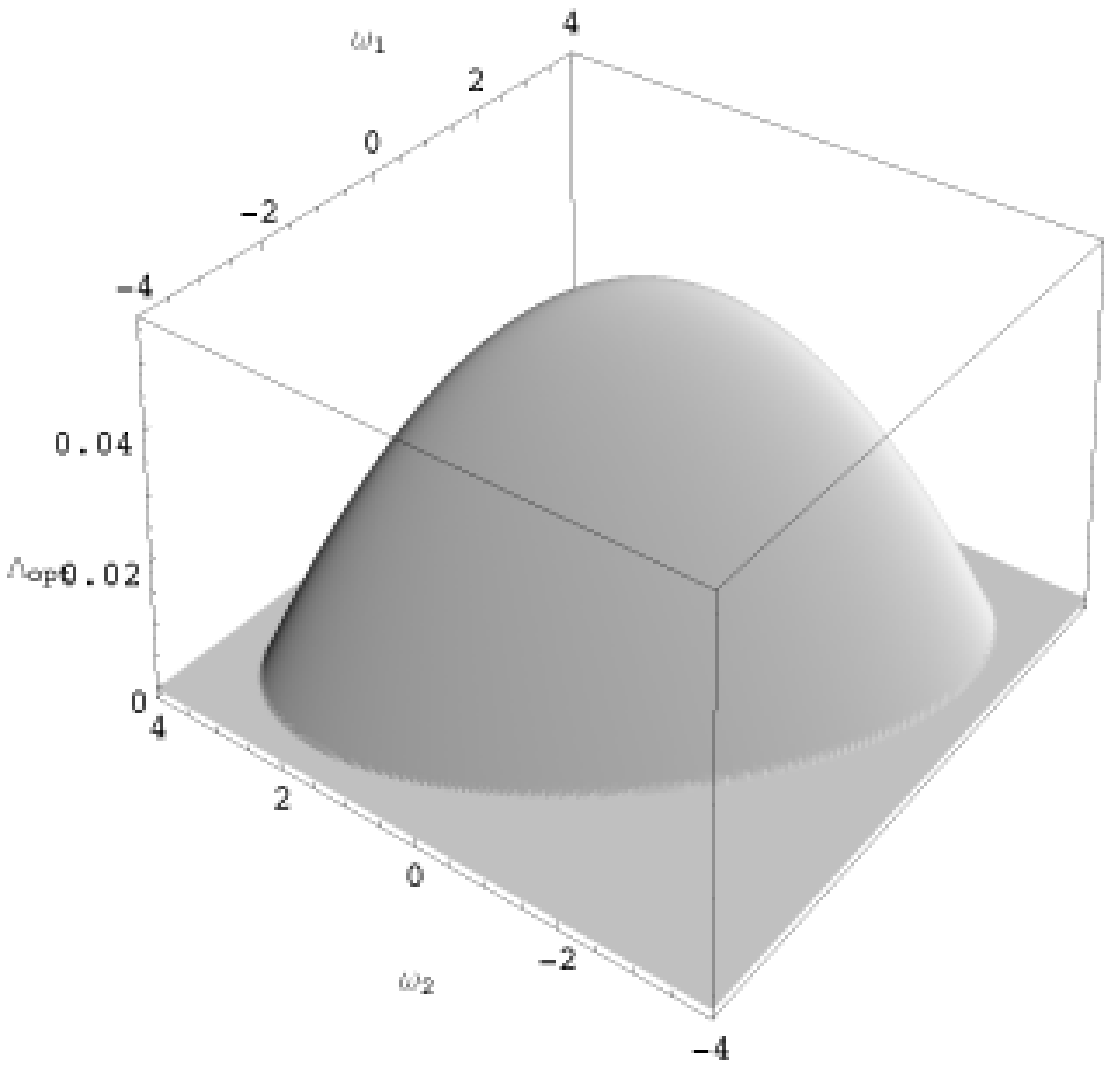}\hspace{1cm}\includegraphics[width=2.5in]{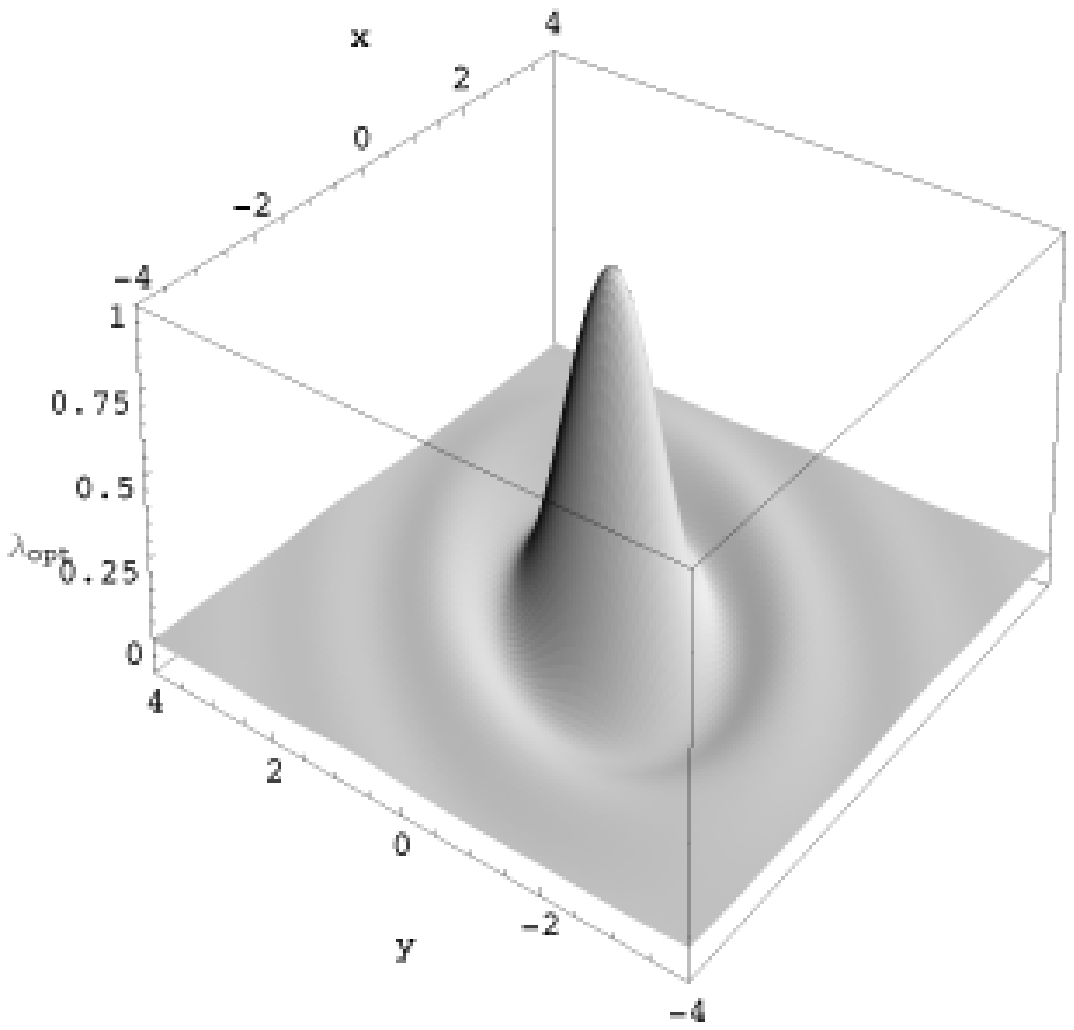}}
\caption{Plots of $\Lambda_{\text{opt}}(\omega_1,\omega_2)$ and
$\lambda_{\text{opt}}(x,y)$.}
\end{figure}

The kernel $\Lambda_{\text{opt}}(\omega_1,\omega_2)$ can be derived
from \eqref{eqn:fsym} (unlike \eqref{eqn:decomp} as proved by
Theorem 1) by symmetrizing $f(x,y)=1-x^2-y^2$ with the symmetric
function $h$ being the arithmetic mean. Since the symmetrized $f$
will be a quadratic form satisfying the appropriate symmetries, it
must be of the form \eqref{eqn:poly} for some $\alpha$ and $\beta$.

Now let $f(x,y)=(1-x^2-y^2)^+=\max(1-x^2-y^2,0)$.  We produce three
new lag-windows from the symmetric functions $h=\prod x_i$,
$h=\max(x_i)$, and $h=\min(x_i)$.

\begin{figure}[H]
\centerline{\includegraphics[width=1.75in]{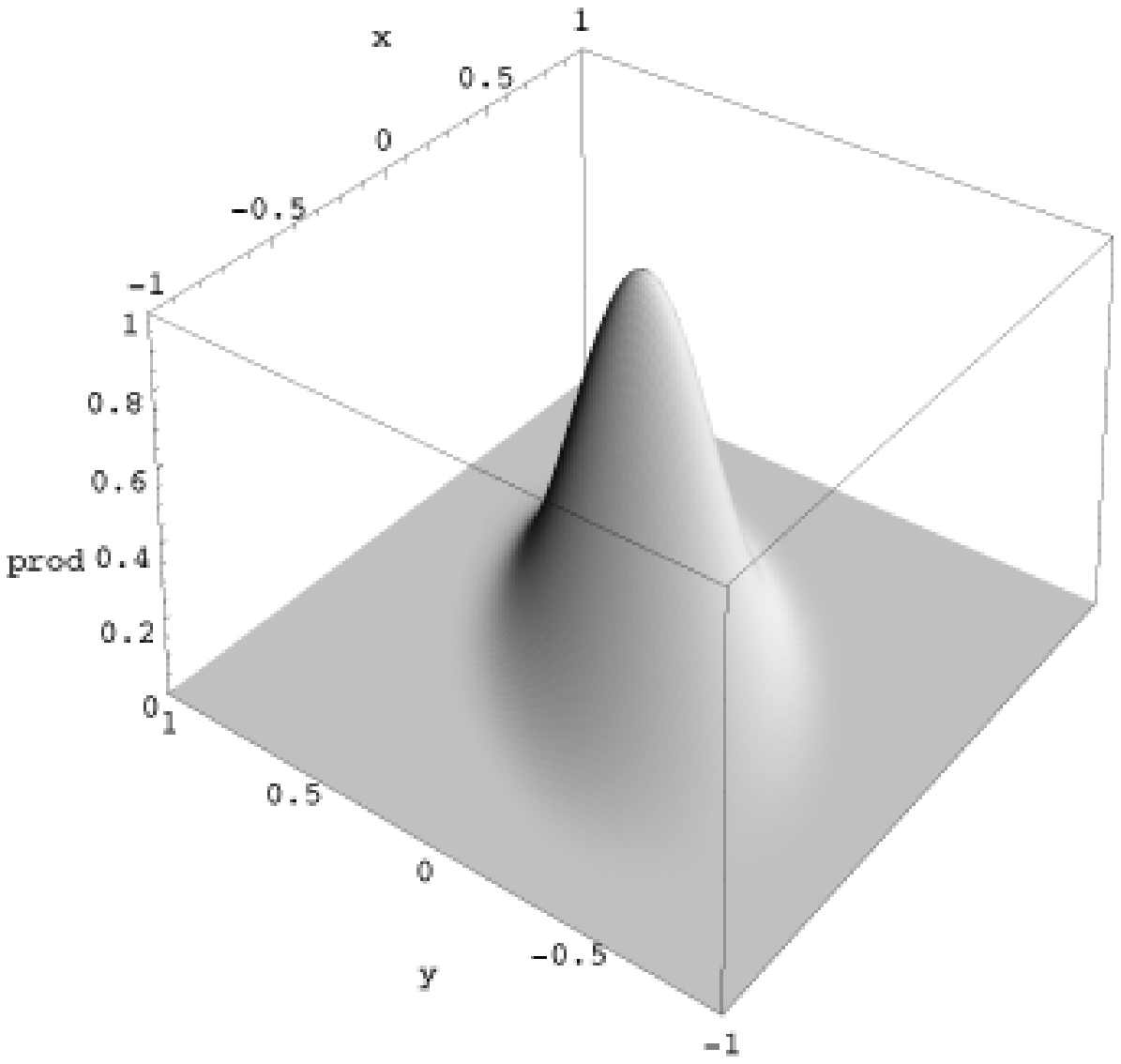}\includegraphics[width=1.75in]{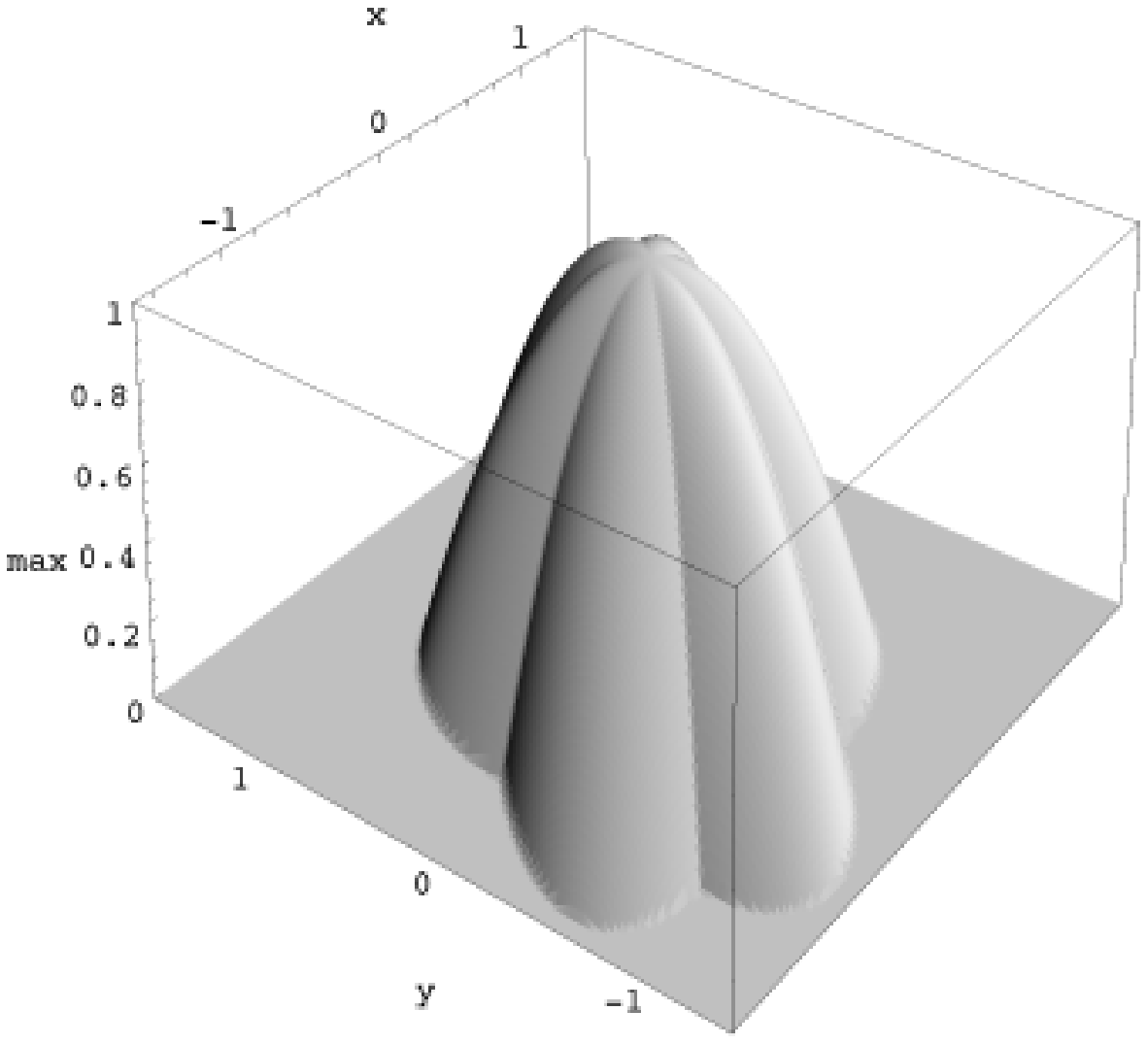}
\includegraphics[width=1.75in]{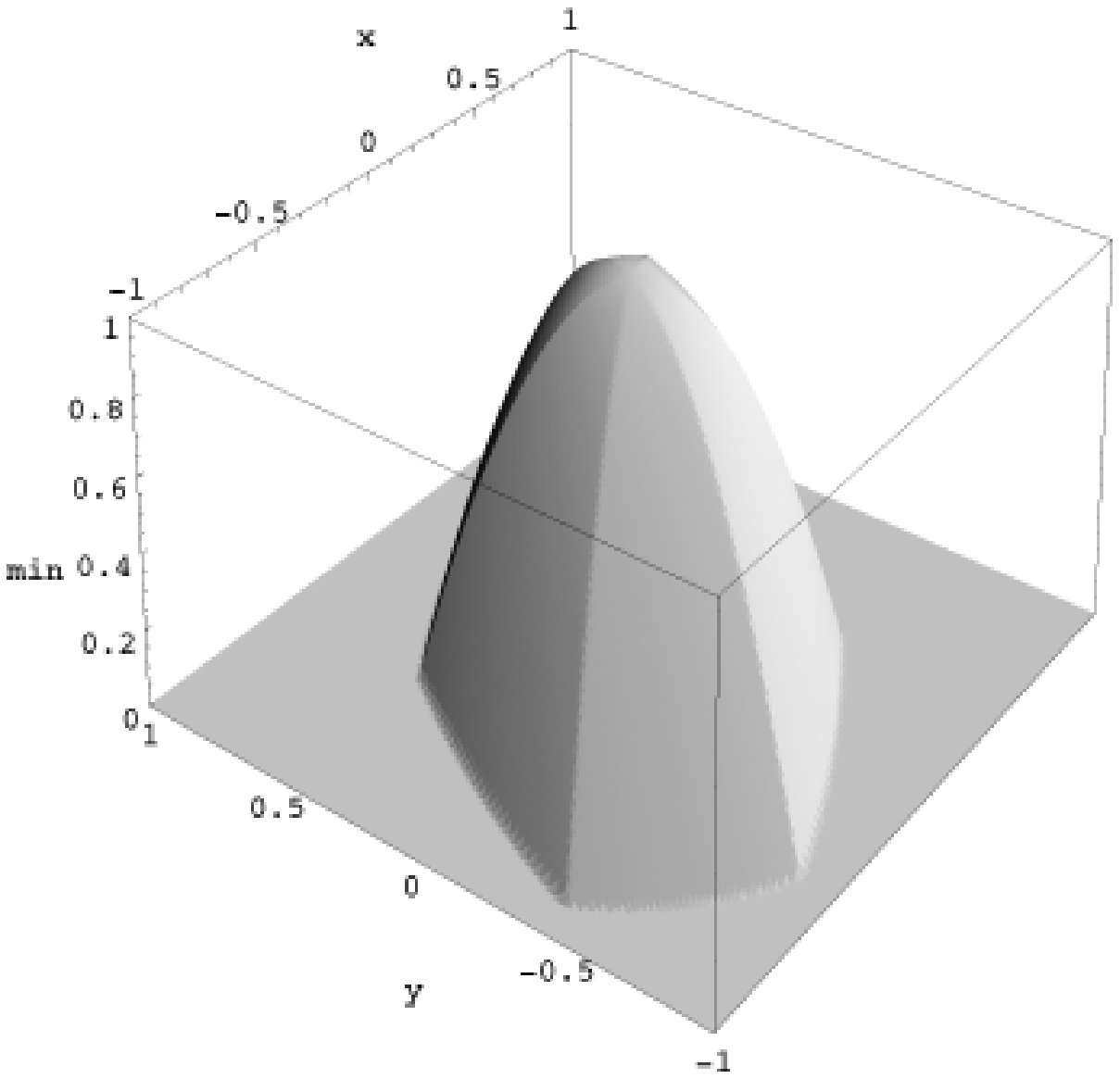}} \caption{$\tilde{f}$ with $h=\prod x_i$,
$h=\max(x_i)$, and $h=\min(x_i)$.}
\end{figure}

We conclude with two flat-top (infinite-order) lag-windows since
estimation with these windows are asymptotically superior
\cite{Berg06}.  The first example is a right pyramidal frustum with
the hexagonal base $|x|+|y|+|x-y|=2$. We let $c\in(0,1)$ be a
scaling parameter that dictates when the frustum becomes flat, that
is, the flat-top boundary is given by $|x|+|y|+|x-y|=2c$.  The
equation of this lag-window is given by
\[
\lambda_{\text{rpf}}(\tau_1,\tau_2)=\frac{1}{1-c}
\lambda_{\text{rp}}(\tau_1,\tau_2)-\frac{c}{1-c}
\lambda_{\text{rp}}\left(\frac{\tau_1}{c},\frac{\tau_2}{c}\right)
\]
where $\lambda_{\text{p}}$ is the equation of the right pyramid with
base $|x|+|y|+|x-y|=2$, i.e.,
\[
\lambda_{\text{rp}}(x,y)=
\begin{cases}
(1-\max(|x|,|y|))^+  ,&-1\le x,y\le 0\ \text{or}\ 0\le x,y\le 1\\
(1-\max(|x+y|,|x-y|))^+, &\text{otherwise}
\end{cases}
\]
The second example is the right conical frustum with elliptical base
$x^2-xy+y^2=1$. Again $c\in(0,1)$ is a scaling parameter, and the
lag-window becomes flat in the ellipse $x^2-xy+y^2=c^2$.  The
equation of this lag-window is given by
\[
\lambda_{\text{rcf}}(\tau_1,\tau_2)=\frac{1}{1-c}
\lambda_{\text{rc}}(\tau_1,\tau_2)-\frac{c}{1-c}
\lambda_{\text{rc}}\left(\frac{\tau_1}{c},\frac{\tau_2}{c}\right)
\]
where $\lambda_{\text{rc}}$ is the equation of the right cone with
base $x^2-xy+y^2=1$, i.e.,
\[
\lambda_{\text{rc}}(x,y)= (1-\sqrt{x^2-xy+y^2})^+
\]
These two examples are plotted below with $c=1/2$ in each case.
\begin{figure}[H]
\centerline{\includegraphics[width=2.5in]{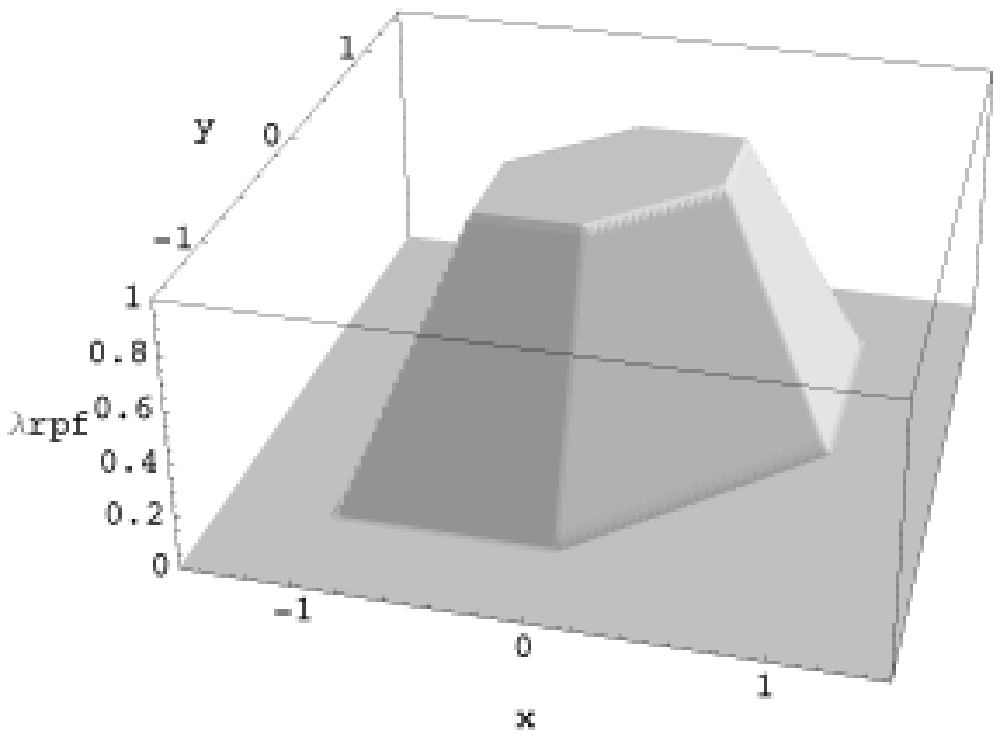}\includegraphics[width=2.5in]{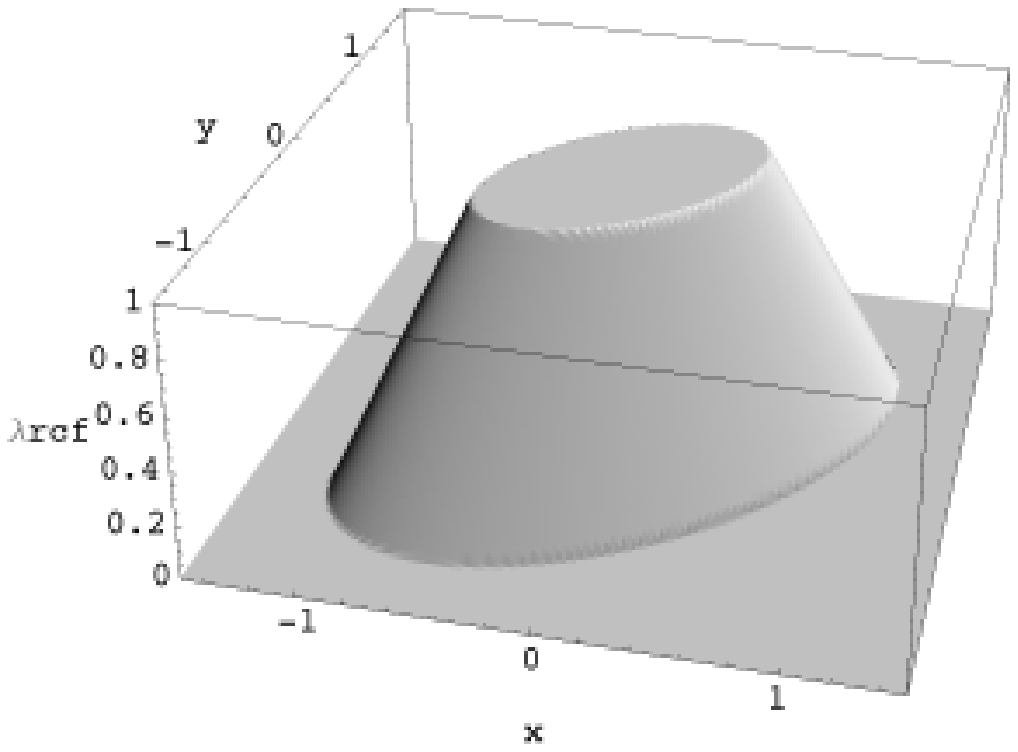}}
\caption{Lag-windows $\lambda_{\text{rpf}}$ and
$\lambda_{\text{rpf}}$.}
\end{figure}

\section{Conclusions}

This paper elicits a deep connection across permutations of $S_k$,
symmetries of lag-windows and kernels, and invertible matrices over
$\bR$.  The structure provided in group theory gives a clear
understanding to many of the basic elements of higher-order spectral
analysis.  In particular, the symmetries of the auto-cumulant
function are now well understood, a general prescription for
producing multivariate lag-windows and kernels is given, and the
Gabr-Rao optimal kernel is easily generalized with the aid of basic
properties of the symmetric group.

%





\nocite{Lu96,nikias93,Nikias87,Sasaki75,Turkbeyler00,Yao01}

\vspace{.5cm}

\noindent {\bf{Acknowledgements.}}  The author thanks John Hall for
his helpful formulation of the representation produced in
(\ref{eqn:compact}) and also thanks Jason Bandlow, Gregg Musiker,
and Oded Yacobi for their helpful insights and discussions on the
symmetric group and group representations.


\bibliography{symmetry}

\end{document}